\documentclass[12pt]{amsart}
\usepackage{amscd}
\def\NZQ{\Bbb}               

\def\QQ{{\NZQ Q}}
\def\ZZ{{\NZQ Z}}

\def\FF{{\NZQ F}}
\def\GG{{\NZQ G}}
\def\HH{{\NZQ H}}
\def\TT{{\NZQ T}}
\def\DD{{\NZQ D}}
\def\frk{\frak}               

\def\Phi{{\frk n}}
\def\Phi{{\frk N}}
\def\opn#1#2{\def#1{\operatorname{#2}}} 
\opn\chara{char} \opn\length{\ell} \opn\pd{pd} \opn\rk{rk}
\opn\projdim{proj\,dim} \opn\injdim{inj\,dim} \opn\rank{rank}
\opn\depth{depth} \opn\grade{grade} \opn\height{height}
\opn\embdim{emb\,dim} \opn\codim{codim}

\opn\Tr{Tr} \opn\bigrank{big\,rank}
\opn\superheight{superheight}\opn\lcm{lcm}
\opn\trdeg{tr\,deg}
\opn\reg{reg} \opn\lreg{lreg} \opn\ini{in} \opn\lpd{lpd}
\opn\size{size} \opn\Pf{Pf} \opn\GL{GL} \opn\SL{SL} \opn\mod{mod}
\opn\ord{ord} \opn\Gin{Gin}
\opn\Hilb{Hilb}\opn\adeg{adeg}\opn\std{std}\opn\ip{infpt}
\opn\div{div} \opn\Div{Div} \opn\cl{cl} \opn\Cl{Cl}
\opn\Spec{Spec} \opn\Supp{Supp} \opn\supp{supp} \opn\Sing{Sing}
\opn\Ass{Ass} \opn\Min{Min}
\opn\Ann{Ann} \opn\Rad{Rad} \opn\Soc{Soc}
\opn\Syz{Syz} \opn\Im{Im} \opn\Ker{Ker} \opn\Coker{Coker}
\opn\Am{Am} \opn\Hom{Hom} \opn\Tor{Tor} \opn\Ext{Ext}
\opn\End{End} \opn\Aut{Aut} \opn\id{id}

\opn\nat{nat}
\opn\pff{pf}
\opn\Pf{Pf} \opn\GL{GL} \opn\SL{SL} \opn\mod{mod} \opn\ord{ord}
\opn\Gin{Gin} \opn\Hilb{Hilb}
\opn\aff{aff} \opn\con{conv} \opn\relint{relint} \opn\st{st}
\opn\lk{lk} \opn\cn{cn} \opn\core{core} \opn\vol{vol}
\opn\link{link} \opn\star{star}
\opn\gr{gr}

\def\pot#1#2{#1[\kern-0.28ex[#2]\kern-0.28ex]}

\opn\dirlim{\underrightarrow{\lim}}
\opn\inivlim{\underleftarrow{\lim}}
\let\sect=\cap
\let\dirsum=\oplus
\let\tensor=\otimes
\let\iso=\cong
\let\Union=\bigcup

\let\Dirsum=\bigoplus

\let\to=\rightarrow
\let\To=\longrightarrow
\def\Implies{\ifmmode\Longrightarrow \else
        \unskip${}\Longrightarrow{}$\ignorespaces\fi}
\def\implies{\ifmmode\Rightarrow \else
        \unskip${}\Rightarrow{}$\ignorespaces\fi}
\def\iff{\ifmmode\Longleftrightarrow \else
        \unskip${}\Longleftrightarrow{}$\ignorespaces\fi}

\let\:=\colon
\newtheorem{Theorem}{Theorem}[section]

\newtheorem{Corollary}[Theorem]{Corollary}
\newtheorem{Proposition}[Theorem]{Proposition}

\let\epsilon\varepsilon
\let\phi=\varphi
\let\kappa=\varkappa
\def\qed{\ifhmode\textqed\fi
      \ifmmode\ifinner\quad\qedsymbol\else\dispqed\fi\fi}
\def\textqed{\unskip\nobreak\penalty50
       \hskip2em\hbox{}\nobreak\hfil\qedsymbol
       \parfillskip=0pt \finalhyphendemerits=0}
\def\dispqed{\rlap{\qquad\qedsymbol}}

\opn\dis{dis}
\def\pnt{{\raise0.5mm\hbox{\large\bf.}}}

\opn\Lex{Lex}

\begin{document}

\title{A generalization of the Taylor complex construction}\maketitle
\begin{center}{\sc J\"urgen Herzog}
\end{center}
\begin{center}   Fachbereich Mathematik und Informatik \\
  Universit\"at Duisburg-Essen, Campus Essen\\
   45117 Essen, Germany\\
e-mail: juergen.herzog@uni-essen.de
\end{center}
\begin{abstract}
Given  multigraded free resolutions of two monomial ideals we
construct a multigraded free resolution of the sum of the two
ideals.
\end{abstract}

\section*{Introduction}
Let $K$ be a field, $S=K[x_1,\ldots, x_n]$ the polynomial ring in
$n$ variables over $K$,   and let $I$ and $J$ be two monomial
ideals in $S$. Suppose that $\FF$ is a multigraded  free
$S$-resolution of $S/I$ and $\GG$ a multigraded  free
$S$-resolution of $S/J$. In this note we construct a  multigraded
free resolution of $S/(I+J)$ which we denote by $\FF*\GG$. It
follows from our construction that $\beta_i(S/(I+J))\leq
\sum_{j=0}^i\beta_{j}(S/I)\beta_{i-j}(S/J)$ for all $i\geq 0$.
Here $\beta_i(M)$ denotes the $i$th Betti number of a graded
$S$-module $M$, that is, the $K$-dimension of $\Tor_i^S(K,M)$.

The inequality for the Betti-numbers  implies in particular that
$\projdim(I+J)\leq \projdim(I)+\projdim(J)+1$. The numerical data
of the complex $\FF*\GG$  also yield the inequality $\reg(I+J)\leq
\reg(I)+\reg(J)-1$. Similar inequalities hold for the projective
dimension and the regularity of $I\sect J$, see Section 3. The
inequality for the regularity has first been conjectured  by
Terai \cite{T}. He also proved this inequality in a special case.

In the squarefree case these inequalities  have first been proved
by Kalai and Meshulam \cite{KM}. The construction of the complex
$\FF*\GG$ was inspired by the work of Kalai and Meshulam. In fact,
the first author informed me that the above mentioned inequalities
for the projective dimension and the regularity of sums and
intersections of squarefree monomial ideals follow from certain
inequalities proved in \cite{KM} concerning the $d$-Leray
properties of the union and intersection of simplicial complexes.
Thus our construction provides an algebraic explanation of   these
inequalities.

One should note that for example the inequality $\projdim(I+J)\leq
\projdim(I)+\projdim(J)+1$, as well as all the other inequalities,
are  wrong for arbitrary graded ideals.

We would also like to mention that the Taylor resolution
(cf.\cite{Ei}) is a special case of our construction. The Taylor
resolution  is a multigraded free resolution for monomial ideals.
It has a uniform structure, but in most cases,  the Taylor
resolution is non-minimal.  In the frame of our construction the
Taylor resolution  can be described as follows: if $I\subset S$ is
a monomial ideal with the minimal set of monomial generators
$\{u_1,\ldots, u_r\}$, and $\FF_j$ is the graded minimal free
resolution of the principal ideal $(u_j)$ for $j=1,\ldots,r$, then
$\FF_1*\FF_2*\cdots*\FF_r$ is the Taylor resolution of $S/I$.

\section{The construction}
Let $K$ be field, $S=K[x_1,\ldots, x_n]$ a polynomial ring and
$I\subset S$ a monomial ideal. Then $S/I$ admits a multigraded
minimal free $S$-resolution
\[
\begin{CD}
\FF\: 0@>>> F_p@>\phi_p >> F_{p-1}@> \phi_{p-1} >> \cdots @>\phi_2
>> F_1@>\phi_1>> F_0 @>>> 0,
\end{CD}
\]
that is, one has
\begin{enumerate}
\item[(i)] $H_0(\FF)=S/I$; \item[(ii)] $F_i=\Dirsum_jS(-a_{ij})$
with $a_{ij}\in  \ZZ^n$ for all $i$; \item[(iii)] the
differentials $\phi_i$ are homomorphisms of multigraded modules.
\end{enumerate}

We define a partial order on $\ZZ^n$ by saying that $b\leq a$ for
$a,b\in\ZZ^n$, if $b$ is componentwise less than $a$.

For all $i$ let $B_i$ be a multihomogeneous basis of $F_i$. Then
$F_i=\Dirsum_{g\in B_i}Sg$, and the differential $\phi_i\: F_i\to
F_{i-1}$ can be described by the equations
\begin{eqnarray}
\phi_i(g)=\sum_{h\in B_{i-1}}a_{gh}u_{gh}h,
\end{eqnarray}
where $a_{gh}\in K$ with $a_{gh}=0$ if $\deg g< \deg h$, and where
the coefficient $u_{gh}$ is the unique monomial in $S$ with $\deg
g=\deg u_{gh}+\deg h$ whenever  $\deg g\geq \deg h$.

For a homogeneous element $f$ in a multigraded module we denote by
$u_f$ the unique monomial with $\deg u_f=\deg f$. Then for all
$i$, all $g\in B_i$  and  all $h\in B_{i-1}$ with $a_{gh}\neq 0$
we have $u_{gh}=u_g/u_h$.

Now let $J\subset S$ be another  monomial ideal with minimal
multigraded free resolution $\GG$ whose differential $\psi$ is
given by
\begin{eqnarray}
\psi_i(e)=\sum_{f\in B'_{i-1}}b_{ef}u_{ef}f\quad\text{for}\quad
e\in B_i',
\end{eqnarray}
where for each $i$, $B'_i$ is a homogeneous basis of $G_i$.

We are going to construct an acyclic multigraded complex $\FF*\GG$
of free $S$-modules which provides a free resolution of $S/(I+J)$.

Let $F$ be a multigraded free $S$-module with homogeneous basis
$B$ and $G$ a multigraded free $S$-module with homogeneous basis
$B'$. We let $F*G$ be the multigraded free $S$-module with a basis
given by the symbols $f*g$ where  $f\in B$ and $g\in B'$.  The
multidegree of $f*g$ is defined to be
\[
\deg f*g =[u_f,u_g],
\]
where for two monomials $u,v\in S$ we denote by $[u,v]$ the least
common multiple of $u$ and $v$. Denoting  by $(u,v)$ the greatest
common divisor of $u$ and $v$, the map
\begin{eqnarray}
\label{mono} j\: F\tensor G\to F*G, \quad f\tensor g\mapsto
(u_f,u_g) f*g
\end{eqnarray}
is a multigraded monomorphism.

\medskip
Now we are ready to define $\FF*\GG$: we let
\[
(\FF*\GG)_i=\Dirsum_{j,k\atop j+k=i}F_j*G_k,
\]
and define the differential
\[
\alpha_i\: (\FF*\GG)_i\To (\FF*\GG)_{i-1}
\]
by the equation
\[
\alpha_i(g*e)=\sum_{h\in B_{j-1}}a_{gh}u_{ghe}h*e+(-1)^j\sum_{f\in
B'_{k-1}}b_{ef}u_{gef}g*f
\]
where  $g\in B_j$ and $e\in B'_k$ with $j+k=i$. Here
\[
u_{ghe}=[u_g,u_e]/[u_h,u_e]\quad \text{and}\quad
u_{gef}=[u_g,u_e]/[u_g,u_f].
\]
Extending the multigraded monomorphism (\ref{mono}) naturally to
$\FF\tensor \GG$ we obtain a  monomorphism of multigraded modules
\begin{eqnarray}
\label{comparison} j\: \FF\tensor \GG \To \FF*\GG
\end{eqnarray}
with the property that $\alpha\circ j=j\circ\partial$, where
$\partial$ denotes the differential of $\FF\tensor \GG$. Since $j$
becomes an isomorphism after localization with respect to all
variables,  it follows that $\alpha\circ \alpha=0$, so that
$\FF*\GG$ is a complex of multigraded $S$-modules.

\section{Acyclicity}
The aim  of this section is to prove the following
\begin{Theorem}
\label{main} Let $I$ and $J$ be  monomial ideals in
$S=K[x_1,\ldots,x_n]$ with  multigraded free $S$-resolutions $\FF$
and $\GG$, respectively. Then $\FF*\GG$ is a multigraded  free
$S$-resolution of $S/(I+J)$.
\end{Theorem}

For the proof of this theorem we need to consider polarization of
monomial ideals. Let $I=(u_1,\ldots, u_m)$ with
$u_j=x_1^{a_{j1}}\cdots x_n^{a_{jn}}$, and let $a$ be the maximum
of the exponents $a_{ij}$. We denote by $\sup u_j$ the set of
elements $i\in [n]$ with $a_{ji}\neq 0$. Consider the polynomial
ring $T$ over $S$ in the variables $y_{ik}$, $i=1,\ldots,n$ and
$k=1,\ldots,r$ with $r\geq a$. The {\em polarization of $I$}  is
the squarefree monomial ideal $I^\wp\subset T$ whose generators
are the monomials
\[
u_j^\wp=\prod_{i\in\sup u_j}y_{i1}y_{i2}\cdots y_{i,a_{ji}},\quad
j=1,\ldots,m.
\]
It is known that
\begin{eqnarray}
\label{regular} \qquad {\bf z}= x_1-y_{11},\ldots,
x_1-y_{1r},x_2-y_{21},\ldots, x_2-y_{2r}, \ldots,
x_n-y_{n1},\ldots x_n-y_{nr}
\end{eqnarray}
is a regular sequence  on $T/I^\wp$ with $(T/I^\wp)/({\bf
z})(T/I^\wp)\iso S/I$.

\medskip
Let $\FF$ be a minimal multigraded free resolution of $S/I$. We
shall need the following result of Sbarra \cite{Sb}, whose proof
we indicate for the convenience of the reader.

\begin{Proposition}
\label{polarize} Let $\FF^\wp$ be a minimal multigraded free
resolution of $T/I^\wp$ and for each $i$, let $B^\wp_i$ be a
multihomogeneous basis of $\FF_i^\wp$. Then for each $i$, there
exists a multihomogeneous basis $B_i$ of $F_i$ and a bijection
$B^\wp_i\to B_i$, $f\mapsto \bar{f}$ with the property that
$u_f=u_{\bar f}^\wp$ for all $f\in B_i^\wp$. In other words, the
shifts in the resolution of $\FF^\wp$ are obtained from the shifts
in $\FF$ by polarization.
\end{Proposition}

\begin{proof}
Notice that $\FF ^\wp/({\bf z})\FF^\wp$ is a minimal graded free
$S$-resolution of $S/I$ since the  sequence ${\bf z}$ (see
(\ref{regular})) is regular on $T/I^{\wp}$. With respect to the
coarse multigrading on $T$  which  assigns to each $y_{ik}$ and to
each $x_i$ the multidegree $\epsilon_i$ where $\epsilon_i$ is
$i$th vector of the canonical basis of $\QQ^n$, the sequence $\bf
z$ is even homogeneous, so that $\FF ^\wp/({\bf z})\FF^\wp$ is a
multigraded complex of $S$-modules, and hence  as a multigraded
complex is isomorphic to $\FF$. Thus we may identify $\FF
^\wp/({\bf z})\FF^\wp$ with $\FF$.

Let  $f\in B_i^\wp$. We denote the residue class of $f$ in  $\FF
^\wp/({\bf z})\FF^\wp$ by $\bar{f}$, and set $B_i=\{\bar{f}\: f\in
B_i^\wp\}$. Then for all $i\geq 0$, $B_i$ is a multihomogeneous
basis of $F_i$.

Since $I^\wp$ is a squarefree monomial ideal, each $u_f$ is a
squarefee monomial. In other words, $u_f=\prod_{i=1}^n\prod_{j\in
A_i}y_{ij}$ with certain $A_i\subset \{1,\ldots,r\}$.  It follows
that $u_{\bar{f}}=\prod_{i=1}^nx_i^{|A_i|}$.  If we can show that
each $A_i$ is of the form $A_i=\{1,\ldots,k_i\}$ for some $k_i$,
then $u_{\bar{f}}^\wp=u_f$, as desired.

Since the Taylor complex $\TT$ of $I^\wp$ is a multigraded free
$T$-resolution of $T/I^\wp$, while $\FF^\wp$ is a {\em minimal}
multigraded free $T$-resolution of $T/I^\wp$, we conclude that
$\FF^\wp$ is isomorphic to a multigraded direct summand of $\TT$.
Let $G(I^\wp)$ be the unique minimal monomial set of generators of
$I^\wp$. The shifts of $\TT$ are the least common multiples of
subsets of $G(I^\wp)$. Since each of the generators of $I^\wp$ is
of the form $\prod_{i=1}^n\prod_{j=1}^{k_i}y_{ij}$, it follows
that all shifts of $\TT$, and hence all shifts of $\FF^\wp$ are of
the same form, as desired.
\end{proof}

\begin{proof}[Proof of Theorem \ref{main}] Obviously one has $H_0(\FF*\GG)=S/(I+J)$. In order to show that $\FF*\GG$ is acyclic, we first treat the case that
$I$ and $J$ are squarefree  monomial ideals and that the
resolutions $\FF$ and $\GG$ are minimal. We consider the following
complex filtration of $\FF*\GG$:
\[
0={\mathcal F}^0(\FF*\GG)\subset {\mathcal F}^1(\FF*\GG)\subset
\cdots \subset {\mathcal F}^n(\FF*\GG)=\FF*\GG,
\]
where
\[
{\mathcal F}^j(\FF*\GG)=\Dirsum_{i\leq j}\FF*G_i.
\]
The factor complexes ${\mathcal F}^j(\FF*\GG)/{\mathcal
F}^{j-1}(\FF*\GG)$ are isomorphic to $\FF*G_j$ with differential
given by
\begin{eqnarray}
\label{factor} F_i*G_j\To F_{i-1}*G_j,\quad g*e\mapsto
 \sum_{h\in B_{i-1}}a_{gh}u_{ghe}h*e
\end{eqnarray}
Here we use the assumptions and notation introduced in the
previous section.

The $E^2$-terms of the first quadrant spectral sequence induced by
the filtration ${\mathcal F}$  are given by the homology of the
factor complexes, that is,
\[
E^2_{i,j}=H_j(\FF*G_i) \quad \text{for all}\quad i,j.
\]
We claim that each of  these factor complexes is acyclic. To this
end we first notice that   $\FF*G_j$ is the direct sum of the
complexes $\FF*Se$ with $e\in B'_j$. In other words,
\[
\FF*G_j\iso \Dirsum_{e\in B'_j}\FF*Se.
\]
Thus it suffices to show that each of the complexes $\FF*Se$ is
acyclic.

The complex homomorphism (\ref{comparison}) restricts to  the
complex homomorphism
\begin{eqnarray}
\label{restriction} j\: \FF\To \FF*Se,\quad g\mapsto
(u_g,u_e)g*e\quad \text{for}\quad g\in\Union_iB_i.
\end{eqnarray}
Thus  after localization we have an isomorphism of complexes
\[
(\FF)_{u_e}\iso (\FF*Se)_{u_e}.
\]
In particular, it follows that $(\FF*Se)_{u_e}$ is acyclic.

Without loss of generality we may assume that
$u_e=\prod_{i=k+1}^nx_i$. Now since the differentials of $\FF*Se$
are given by
\[
F_i*Se\to F_{i-1}*Se,\quad g*e\mapsto \sum_{h\in
B_{i-1}}a_{gh}u_{ghe}h*e\quad \text{with}\quad
u_{ghe}=[u_g,u_e]/[u_h,u_e],
\]
we see that all the monomial entries $u_{ghe}$ of the
differentials are monomials in $S'=K[x_1,\ldots,x_k]$, so that
$\FF*Se\iso \HH\tensor_{S'}S$ where $\HH$ is a multigraded complex
of free $S'$-modules. (Here is where we use that $I$ and $J$ are
squarefree). Hence
\[
(\FF*Se)_{u_e}\iso(\HH\tensor_{S'}S)\tensor_SS_{u_e}\iso
\HH\tensor_{S'}S'[x_{k+1}^{\pm 1},\ldots, x_n^{\pm 1}].
\]
Since  $(\FF*Se)_{u_e}$ is acyclic and the since the extension
$S'\subset S'[x_{k+1}^{\pm 1},\ldots, x_n^{\pm 1}]$ is faithfully
flat, it follows that $\HH$ is acyclic. Then, using the fact that
the extension $S'\subset S$ is flat, we conclude that $\FF*Sg\iso
\HH\tensor_{S'}S$ is acyclic, as desired.

Note that $H_0(\FF*Sg)=(S/I_g)g$, where $I_g$ is generated by the
monomials $[u,u_g]/u_g$ with $u\in G(I)$. Here we denote  as usual
by $G(L)$ the unique minimal set of monomial generators of a
monomial ideal $L$. Thus our calculations have shown that
\[
E^2_{i,j}=\left\{ \begin{array}{lll} 0, & \mbox{if} & j>0,\\
\Dirsum_{g\in B'_i}(S/I_g)g, & \mbox{if} & j=0.
\end{array} \right.
\]
Therefore $\FF*\GG$ will be acyclic if the complex
\[
\label{induced}\tilde{\GG}\: 0\To \Dirsum_{g\in B'_q}(S/I_g)g\To
\cdots \To\Dirsum_{g\in B'_2}(S/I_g)g\To \Dirsum_{g\in
B'_1}(S/I_g)g\To S/I\To 0,
\]
is acyclic, where the differentials $\tilde{\psi_i}$ of
$\tilde{\GG}$ are induced by those of $\GG$.

In order to prove the acyclicity of $\tilde{\GG}$, let $i>0$ and
$z\in \tilde{G_i}$ be  a multihomogeneous element with
$\tilde{\psi_i}(z)=0$, and let $w\in G_i$ be a multihomogeneous
element with $\epsilon(w)=z$, where $\epsilon\: G_i\to
\tilde{G_i}$ is the canonical epimorphism. Then
\[
\psi_i(w)\in \Dirsum_{h\in B'_{i-1}}I_hh,
\]
and we have to show that $w\in \Im(\psi_{i+1})+ \Dirsum_{g\in
B'_{i}}I_gg$.

We have  $w=\sum_gc_gw_gg$ with $c_g\in K$ and $w_g$ a monomial
with $\deg w_g+\deg g= \deg w$ for all $g$ with $c_g\neq 0$. Then
\[
\psi_i(w)=\sum_gc_gw_g(\sum_hb_{gh}u_{gh}h)=\sum_h(\sum_gc_gb_{gh}w_gu_{gh})h.
\]
The monomial $w_gu_{gh}$  only depends on $w$ and on $h$ (and not
on $g$). We therefore denote it by $v_h$, and obtain
\[
\psi_i(w)=\sum_h(\sum_gc_gb_{gh})v_hh.
\]
If $\psi_i(w)=0$, then $w\in \Im(\psi_{i+1})$ since $\GG$ is
acyclic. Otherwise, there exists $h\in B'_{i-1}$ such that
$\sum_gc_gb_{gh}\neq 0$. For this $h$ one has $v_h\in I_h$, and
there exists $e\in B_i'$ with $c_eb_{eh}\neq 0$. Since
$v_h=w_eu_{eh}=w_e(u_e/u_h)$ it follows that $w_eu_e\in
I_hu_h\subset I$, so that $w_e\in I_e$. Since on the other hand
$\psi_i(z_ee)\in \Dirsum_{h\in B'_{i-1}}I_hh$ for any monomial
$z_e\in I_e$, it follows for $w'=\sum_{g,\; g\neq e}c_gw_gg$ that
\[
\psi_i(w')=\psi(w)-\psi(c_ew_ee)\in \Dirsum_{h\in B'_{i-1}}I_hh.
\]
Hence using induction on the number of summands in $w$ we may
assume that $w'\in \Im(\psi_{i+1})+ \Dirsum_{g\in B'_{i}}I_gg$,
which yields the desired conclusion since $c_ew_ee\in\Dirsum_{g\in
B'_{i}}I_gg$.

\medskip
Next  we consider the case  that $I$ and $J$ are arbitrary
monomial ideals in $S$ and that the resolutions $\FF$ and $\GG$
are minimal. We use polarization, to reduce this more general case
to the case of squarefree monomial ideals.

Assume that  the differential of $\FF^\wp$ is given by
\[
\phi_i(g)=\sum_{h\in B_{i-1}^\wp}a_{gh}u_{gh}h,
\]
and that of $\GG^\wp$ is given by
\[
\psi_i(e)=\sum_{f\in B'^\wp_{i-1}}b_{ef}u_{ef}f.
\]
with multihomogeneous bases $B_{i}^\wp$ and $B'^\wp_{i}$.

 Then the  differential of $\FF^\wp*\GG^\wp$ is given
by
\[
\alpha_i(g*e)=\sum_{h\in B^\wp_{i-1}}a_{gh}u_{g h
e}h*e+(-1)^j\sum_{f\in B'^\wp_{i-1}}b_{ef}u_{g e f}g*f
\]
for $g\in B^\wp_j$ and $f\in B'^\wp_k$, where
\[
u_{g h e}=[u_{g},u_{e}]/[u_{h},u_{e}]\quad \text{and}\quad u_{g e
f}=[u_{g},u_{e}]/[u_{g},u_{f}].
\]

It follows that the differential of $(\FF^\wp*\GG^\wp)/({\bf
z})(\FF^\wp*\GG^\wp)$ is given by
\[
\overline{\alpha}_i(\overline{g*e})=\sum_{h \in
B^\wp_{i-1}}a_{gh}\overline{u}_{g h
e}\overline{h*e}+(-1)^j\sum_{f\in
B'^\wp_{i-1}}b_{ef}\overline{u}_{g e f}\overline{g*f}.
\]
Here $\; \bar{}\; $ denotes  the residue class of an element
modulo ${\bf z}$. Now observe that
\[
\overline{u}_{g h
e}=\overline{[u_{g},u_{e}]/[u_{h},u_{e}]}=\overline{[u_{g},u_{e}]}/\overline{[u_{h},u_{e}]}=[\overline{u}_{g},\overline{u}_{e}]/[\overline{u}_{h},\overline{u}_{e}]
=[u_{\bar{g}},u_{\bar{e}}]/[u_{\bar{h}},u_{\bar{e}}]=u_{\bar{g}\bar{h}\bar{e}},
\]
and similarly, $\overline{u}_{g e f}=u_{\bar{g}\bar{e}\bar{f}}$.
For the third equation we used that
$\overline{[u,v]}=[\bar{u},\bar{v}]$ for  monomials  $u$ and $v$
which are of the form $\prod_{j=1}^n\prod_{i=1}^{k_i} y_{ji}$. But
this condition is satisfied in our case since by Proposition
\ref{polarize} the monomials $u_g$, $u_h$ and $u_e$ are the
polarizations of the monomials $u_{\bar{g}}$, $u_{\bar{h}}$ and
$u_{\bar{e}}$, respectively. For the fourth equation we used, also
shown in Proposition \ref{polarize}, that
$\overline{u}_f=\overline{u^{\wp}_{\bar{f}}}=u_{\bar{f}}$.

The above calculations show  that
\[
(\FF^\wp*\GG^\wp)/({\bf z})(\FF^\wp*\GG^\wp)\To \FF*\GG,\quad
\overline{g*e}\mapsto \bar{g}*\bar{e},
\]
establishes an isomorphism of complexes, as desired.

In the final step of the proof we assume that $I$ and $J$ are
arbitrary monomial ideals but  $\FF$ and $\GG$ are not necessarily
minimal. Then $\FF$ can be written as a direct sum
$\FF=\FF'\dirsum
 \HH$ of multigraded complexes with $\FF'$  a minimal free
 resolution of $S/I$ and  $\HH$ exact. Note that $\HH$ is a
 direct sum of complexes of the form  $\DD: 0\to Sg\to Sh\to 0$
 whose differential maps $g$ to $h$. Since $\DD*\GG$ is isomorphic to the
 mapping cone of the identity on $\GG$ we see that $\DD*\GG$ is
 exact. It follows that $\HH*\GG$ is exact, and consequently
 $\FF*\GG$ has the same homology as $\FF'*\GG$. By the same
 argument we may replace $\GG$ by a minimal multigraded free
 resolution of $\GG'$ and thus obtain that
\[
H_i(\FF*\GG)\iso H_i(\FF'*\GG')=\left\{ \begin{array}{lll} 0, & \mbox{if} & i>0,\\
S/(I+J), & \mbox{if} & i=0.
\end{array} \right.
\]
\end{proof}

\section{Consequences}

The product $\FF*\GG$ defined for multigraded free resolutions is
associative, that is, we have
\[
(\FF*\GG)*\HH\iso \FF*(\GG*\HH)
\]
for any three multigraded free resolutions. Thus if $I_j\subset S$
is a monomial ideal and   $\FF_j$ is a  multigraded free
$S$-resolution of $S/I_j$ for $j=1,\ldots,r$, then
\[
\FF_1*\FF_2*\cdots * \FF_r
\]
is multigraded free $S$-resolution of $S/(I_1+I_2+\cdots +I_r)$.

Consider the following special case: let $I$ be a monomial ideal
with unique minimal monomial set of generators $G(I)=\{u_1,\cdots,
u_r\}$, and set $I_j=(u_j)$ and  $a_j=\deg u_j$ for
$j=1,\ldots,r$. Then
\[
\begin{CD}
\FF_j\: 0@>>> S(-a_j) @> u_j >> S@>>> 0
\end{CD}
\]
is a multigraded free $S$-resolution of $S/I_j$, and so $\TT =
\FF_1*\FF_2*\cdots
* \FF_r$ is a multigraded free $S$-resolution of $S/I$.
Indeed, $\TT$ is the well-known Taylor resolution of $S/I$ (cf.\
\cite[Exercise 17.11]{Ei}).

An obvious consequence of our construction is the following

\begin{Corollary}
\label{bound} Let $I$ and $J$ be monomial ideals in $S$. Then
\[
\beta_i(S/(I+J))\leq \sum_{j=0}^i\beta_{j}(S/I)\beta_{i-j}(S/J).
\]
\end{Corollary}

For a graded ideal $L\subset S$  we set $M_i(L)=\max\{j\:
\Tor_i^S(K,L)_j\neq 0\}$. In other words, $M_i(L)$ is the highest
shift in the $i$th step of the graded minimal free  resolution of
$L$. Furthermore we set $\reg_i(L)=M_i(L)-i$ for $i\geq 0$ and
$\reg(L)_{-1}=0$. The {\em regularity} of $L$ is then defined to
be $\max\{\reg_i(L)\: i\geq 0\}$.

Corollary \ref{bound} implies the inequalities (a) described in
the next corollaries. Inequality \ref{conj}(b) concerning the
regularity was conjectured by Terai \cite{T} and proved in a
special case.

\begin{Corollary}[Kalai, Meshulam]
\label{conj} Let $I$ and $J$ be monomial ideals. Then
\begin{enumerate}
\item[(a)]  $\projdim(I+J)\leq \projdim(I)+\projdim(J)+1$;
\item[(b)]  $\reg(I+J)\leq \reg(I)+\reg(J)-1$.
\end{enumerate}
\end{Corollary}

\begin{proof} It remains to prove statement (b). Since $\FF*\GG$
is a possibly non-minimal graded free resolution of $S/(I+J)$, we
see that $M_i(I+J)$ is less than or equal to the maximal
$\ZZ$-degree of a generator of $(\FF*\GG)_{i+1}$.

Since $(\FF*\GG)_{i+1}=\Dirsum_{j+k=i+1}F_j*G_k$ and since $\deg
f*g\leq \deg f+\deg g$ for all homogeneous elements $f$ and $g$,
it follows that

\begin{eqnarray*}
\reg_i(I+J)&\leq &\max_{j+k=i+1}\{M_{j-1}(I)+M_{k-1}(J)\}-i\\
 &=& \max_{j+k=i+1}\{\reg_{j-1}(I)+\reg_{k-1}(J)\}-1.
\end{eqnarray*}
This implies the desired inequality.
\end{proof}

From the exact sequence
\[
0\To I\sect J\To I\dirsum J\To I+J\To 0
\]
and Corollary \ref{conj} we deduce the following inequalities

\begin{Corollary}[Kalai, Meshulam]
\label{intersect} Let $I$ and $J$ be monomial ideals. Then
\begin{enumerate}
\item[(a)] $\projdim(I\sect J)\leq \projdim(I)+\projdim(J)$;

\item[(b)] $\reg(I\sect J)\leq \reg(I)+\reg(J)$.
\end{enumerate}
\end{Corollary}

In the case of monomial complete intersections, a proof of
inequality (b) is given  by Marc Chardin, Nguyen Cong Minh and Ngo
Viet Trung \cite{CMT}.

\medskip

For simplicial complexes the inequality of Corollary \ref{bound}
has the following interpretation. Let $\Sigma$ be a simplicial
complex and $U$ a subset of the vertex set. We denote by
$\Sigma_U$ the restriction of $\Sigma$ to $U$, that is, the
simplicial complex with faces $F\in\Sigma$ such that $F\subset U$.
We fix a field and  denote by $\tilde{H}_i(\Sigma)$ the $i$th
reduced simplicial homology of $\Sigma$ with respect to $K$, and
by $\tilde{h}_i(\Sigma)$ the $K$-dimension of
$\tilde{H}_i(\Sigma)$. With this notation we have

\begin{Corollary}
\label{simplicial} Let $\Delta$ and $\Gamma$ be simplicial
complexes on the vertex set $[n]$. Then
\[
\sum_W \tilde{h}_{|W|-i-1}((\Delta\cap \Gamma)_W)\leq
\sum_{j+k=i}(\sum_U \tilde{h}_{|U|-j-1}(\Delta_U))(\sum_ V
\tilde{h}_{|V|-k-1}(\Gamma_V)).
\]
Here the sums are taken over all subsets $U,V, W\subset [n]$
\end{Corollary}

\begin{proof}
The inequality is a consequence of Corollary \ref{bound} and
Hochster's formula
\[
\dim_k \Tor_i^S(K;K[\Sigma])=\sum_W \tilde{h}_{|W|-i-1}(\Sigma_W),
\]
see \cite[Theorem 5.5.1]{BH}.
\end{proof}

\end{document}